%% file: CRnekom.tex
\documentclass[10pt,reqno]{amsart}
\usepackage{amssymb,amscd,amsbsy}
\usepackage{amsthm}
\input{CRinput}

\theoremstyle{remark}

\newtheorem*{rem*}{Remark}

\newcommand\fM{\frak M}

\begin{document}

\newcommand{\vse}{\vspace{.2in}}

\title{Functions of perturbed noncommuting self-adjoint operators}

\maketitle
\begin{center}
\Large
Aleksei Aleksandrov$^{\rm a}$, Fedor Nazarov$^{\rm b}$, Vladimir Peller$^{\rm c}$
\end{center}

\begin{center}
\footnotesize
{\it$^{\rm a}$St-Petersburg Branch, Steklov Institute of Mathematics, Fontanka 27, 191023 St-Petersburg, Russia\\
$^{\rm b}$Department of Mathematics, Kent State University, Kent, OH 44242, USA\\
$^{\rm c}$Department of Mathematics, Michigan State University, East Lansing, MI 48824, USA}
\end{center}

\newcommand{\mt}{{\mathcal T}}

\footnotesize

{\bf Abstract.} We consider functions $f(A,B)$ of noncommuting self-adjoint operators $A$ and $B$ that can be defined in terms of double operator integrals. We prove that if $f$ belongs to the Besov class $B_{\be,1}^1(\R^2)$, then we have the following Lipschitz type estimate in the trace norm:
$\|f(A_1,B_1)-f(A_2,B_2)\|_{\bS_1}\le\const(\|A_1-A_2\|_{\bS_1}+\|B_1-B_2\|_{\bS_1})$. However, the condition $f\in B_{\be,1}^1(\R^2)$ does not imply the Lipschitz type estimate in the operator norm.

\medskip

\begin{center}
{\bf\large Fonctions d'op\'erateurs perturb\'es noncommutants auto-adjoints}
\end{center}

\medskip

{\bf R\'esum\'e.} Nous consid\'erons les fonctions $f(A,B)$ d'op\'erateurs auto-adjoints $A$ et $B$ qui ne commutent pas. Telles fonctions peuvent \^etre d\'efinies
en termes d'int\'egrales doubles op\'eratorielles. Pour $f$ dans l'espace de Besov
$B_{\be,1}^1(\R^2)$ nous obtenons l'estimation lipschitzienne en norme trace: $\|f(A_1,B_1)-f(A_2,B_2)\|_{\bS_1}\le\const(\|A_1-A_2\|_{\bS_1}+\|B_1-B_2\|_{\bS_1})$.
D'autre part la condition $f\in B_{\be,1}^1(\R^2)$ n'implique pas l'estimation lipschitzienne en norme op\'eratorielle.

\normalsize

\

\begin{center}
{\bf\large Version fran\c caise abr\'eg\'ee}
\end{center}

\medskip

Il est bien connu (voir \cite{F}) q'une fonction lipschitzienne sur $\R$ ne doit pas \^etre {\it lipschitzienne op\'eratorielle}, c'est-\`a-dire l'in\'egalit\'e
$\|f(A)-f(B)\|\le\const\|A-B\|$ pour les op\'erateurs auto-adjoints $A$ et $B$ peut \^etre fausse. Dans \cite{Pe1} et \cite{Pe2} des conditions n\'ecessaires et des conditions suffisantes sont donn\'ees pour qu'une fonction $f$ soit lipschitzienne
op\'eratorielle. En particulier,  il est d\'emontr\'e dans \cite{Pe1} et \cite{Pe2} que si $f$ appartient \`a l'espace de Besov $B_{\be1}^1(\R)$, alors $f$ est lipschitzienne op\'eratorielle. Il est aussi bien connu qu'une fonction $f$ est lipschitzenne op\'eratorielle si et seulement si la condition $A-B$ appartient \`a la classe $\bS_1$ (classe trace) implique que $f(A)-f(B)\in\bS_1$ et $\|f(A)-f(B)\|_{\bS_1}\le\const\|A-B\|_{\bS_1}$.

D'autre part, il est d\'emontr\'e dans \cite{AP}   que si $f$ est une fonction h\"olderienne d'ordre $\a$, $0<\a<1$, alors
$\|f(A)-f(B)\|\le\const\|A-B\|^\a$ pour les op\'erateurs auto-adjoints $A$ et $B$.

Les r\'esultats ci-dessus ont \'et\'e g\'eneralis\'es dans \cite{APPS}  au cas de fonctions d'op\'erateurs normaux est dans \cite{NP} au cas de fonctions de $n$-uplets d'op\'erateurs auto-adjoints commutants.

Dans cette note nous consid\'erons les fonctions $f(A,B)$ d'op\'erateurs auto-adjoints $A$ et $B$ qui ne doivent pas commuter. Si une fonction $f$ sur $\R^2$ est un multiplicateurs de Schur, on d\'efinit $f(A,B)$ comme l'int\'egrale double op\'eratorielle
$$
f(A,B)=\iint f(x,y)\,d{E_A}(x)\,dE_B(y)
$$
o\`u $E_A$ et $E_B$ sont le mesures spectrales de $A$ et de $B$ (voir \cite{AP} pour l'information sur les multiplicateurs de Schur et sur les int\'egrales doubles op\'eratorielles). 

Nous d\'emontrons que si $f\in B_{\be1}^1(\R)$, alors l'in\'egalit\'e suivante du type lipschitzenne en norme trace est vrai:
$$
\|f(A_1,B_1)-f(A_2,B_2)\|_{\bS_1}
\le\const\|f\|_{B_{\be,1}^1}(\|A_1-A_2\|_{\bS_1}+\|B_1-B_2\|_{\bS_1}),
$$ 
pour tous les op\'erateurs $A_1$, $A_2$, $B_1$ et $B_2$ auto-adjouints tels que
$A_1-A_2\in\bS_1$ et $B_1-B_2\in\bS_1$. Pour d\'emontrer l'in\'egalit\'e si-dessus nous utilisons la repr\'esentation suivante en termes d'int\'egrale triple 
op\'eratorielle:
$$
f(A_1,B)-f(A_2,B)=\int_\R\int_\R\int_\R\frac{f(x_1,y)-f(x_2,y)}{x_1-x_2}
\,dE_{A_1}(x_1)(A_1-A_2)\,dE_{A_2}(x_2)\,dE_B(y).
$$
L'int\'egrale triple op\'eratorielle est bien d\'efini parce que la diff\'erence divis\'ee appartient au produit tensoriel $L^\be(\R)\otimes_{\rm h}\!L^\be(\R)\otimes^{\rm h}\!L^\be(\R)$ (voir la d\'efinition dans la version anglaise). Plus 
pr\'ecis\'ement, si $f$ est une fonction born\'ee sur $\R^2$ dont la transform\'ee de Fourier a un support dans le disque unit\'e, nous obtenons la repr\'esentation tensorielle suivante:
$$
\frac{f(x_1,y)-f(x_2,y)}{x_1-x_2}=
\sum_{j,k\in\Z}\frac{\sin(x_1-j\pi)}{x_1-j\pi}\cdot\frac{\sin(x_2-k\pi)}{x_2-k\pi}
\cdot\frac{f(j\pi,y)-f(k\pi,y)}{j\pi-k\pi}.
$$
En outre, on a
$$
\sum_{j\in\Z}\frac{\sin^2(x_1-j\pi)}{(x_1-j\pi)^2}
=\sum_{k\in\Z}\frac{\sin^2(x_2-k\pi)}{(x_2-k\pi)^2}=1,
\quad x_1~x_2\in\R,
$$
et
$$
\sup_{y\in\R}\left\|\left\{\frac{f(j\pi,y)-f(k\pi,y)}{j\pi-k\pi}
\right\}_{j,k\in\Z}\right\|_\B\le\const\|f\|_{L^\be(\R)},
$$
o\`u $\B$ est l'espace d'op\'erateurs born\'es sur $\ell^2$. Si $j=k$, alors
$(f(j\pi,y)-f(k\pi,y))(j\pi-k\pi)^{-1}
\df\frac{\partial f}{\partial x}(j\pi,y)$.

D'autre part, il se trouve que pour la m\^eme classe de fonctions c'est impossible d'obtenir une in\'egalit\'e du type lipschitzienne dans la norme op\'eratorielle. 
Plus pr\'ecis\'ement, nous pouvons d\'emontrer qu'il n'y a pas du nombre positif $M$ pour lequel
$$
\|f(A_1,B_1)-f(A_2,B_2)\|\le M\|f\|_{L^\be(\R^2)}(\|A_1-A_2\|+\|B_1-B_2\|),
$$ 
chaque fois que $f$ soit une fonction born\'ee sur $\R$ dont la transform\'ee de Fourier a un support dans le disque unit\'e et $A_1$, $A_2$, $B_1$, $B_2$ soient des op\'erateurs auto-adjoints du rang fini. 

En outre, on peut trouver des op\'erateurs auto-adjoints
du rang fini tels que $\|A_1-A_2\|\le2\pi$, $\|B_1-B_2\|\le2\pi$, et une fonction $f$ sur $\R^2$ dont la transform\'ee de Fourier a un support dans le disque unit\'e et telle que $\|f\|_{L^\be(\R)}\le1$ pour lesquels la norme $\|f(A_1,B_1)-f(A_2,B_2)\|$
soit la plus grande possible.

\medskip

\begin{center}
------------------------------
\end{center}

\medskip

\setcounter{section}{0}
\section{\bf Introduction}

\medskip

It was shown in \cite{F} that a Lipschitz function $f$ on $\R$ does not have to be {\it operator Lipschitz}, i.e., it does not have to satisfy the inequality 
$\|f(A)-f(B)\|\le\const\|A-B\|$ for self-adjoint operators $A$ and $B$. In \cite{Pe1} and \cite{Pe2} it was shown that the condition $f\in B^1_{\be,1}(\R)$ is sufficient for $f$ to be operator Lipschitz (see \cite{Pee} and \cite {Pe4} for information about Besov spaces $B_{p,q}^s$). It is also well known  that $f$ is operator Lipschitz if and only if $f(A)-f(B)$ belongs to trace class $\bS_1$ whenever $A-B\in\bS_1$ and $\|f(A)-f(B)\|_{\bS_1}\le\const\|A-B\|_{\bS_1}$. 
It was shown in \cite{AP} that if $f$ is a H\"older function of order $\a$, 
$0<\a<1$, then it is {\it operator H\"older of order} $\a$, i.e., 
$\|f(A)-f(B)\|\le\const\|A-B\|^\a$ for self-adjoint operators $A$ and $B$.

Later the above results were generalized in \cite{APPS} to functions of normal operators and in \cite{NP} to functions of commuting $n$-tuples of self-adjoint operators.

In this paper we are going to consider functions of noncommuting pairs of self-adjoint operators.

Let $A$ and $B$ be self-adjoint operators on Hilbert space and let $E_A$ and $E_B$ be their spectral measures. Suppose that $f$ is a function of two variables that is defined at least on $\s(A)\times\s(B)$. If $f$ is a Schur multiplier with respect to the pair $(E_A,E_B)$, we define the function $f(A,B)$ of $A$ and $B$ by
\bay
\label{fAB}
f(A,B)\df\iint f(x,y)\,dE_A(x)\,dE_B(y)
\ey
(we refer the reader to \cite{AP} for definitions of Schur multipliers and double operator integrals).
Note that this functional calculus $f\mapsto f(A,B)$ is linear, but not multiplicative.

If we consider functions of bounded operators, without loss of generality we may deal with periodic functions with a sufficiently large period. Clearly, we can rescale the problem and assume that our functions are $2\pi$-periodic in each variable.

If $f$ is a trigonometric polynomial of degree $N$, we can represent $f$ in the form
$$
f(x,y)=\sum_{j=-N}^Ne^{{\rm i}jx}\left(\sum_{k=-N}^N\hat f(j,k)e^{{\rm i}ky}\right).
$$
Thus 
$f$ belongs to the projective tensor product $L^\be\hat\otimes L^\be$ and
$$
\|f\|_{L^\be\hat\otimes L^\be}\le\sum_{j=-N}^N\sup_x
\left|\sum_{k=-N}^N\hat f(j,k)e^{{\rm i}ky}\right|
\le(1+2N)\|f\|_{L^\be}
$$
It follows that every periodic function $f$ of class $B_{\be1}^1(\R^2)$ belongs to 
$L^\be\hat\otimes L^\be$ and the operator $f(A,B)$ is well defined by \rf{fAB}.

\

\section{\bf Lipschitz type estimates in the trace norm}

\

In this paper we use triple operator integrals to estimate functions of perturbed noncommuting operators in trace norm. Let $E_1$, $E_2$, and $E_3$ be spectral measures on measurable spaces $(\X_1,\fM_1)$, $(\X_2,\fM_2)$, and $(\X_3,\fM_3)$.
We say that a function $\Psi$ on $\X_1\times\X_2\times\X_3$ belongs to the
{\it Haagerup tensor product of the spaces} $L^\infty(E_j)$, $j=1,2,3$, (notationally, 
$\Psi\in L^\be(E_1)\otimes_{\rm h}\!L^\be(E_2)\otimes_{\rm h}\!L^\be(E_3)\,\,$)
if $\Psi$ admits a representation
\bay
\label{Htp}
\Psi(x_1,x_2,x_3)=\sum_{j,k\ge0}\a_j(x_1)\b_{jk}(x_2)\g_k(x_3);
\ey
here $\{\a_j\}_{j\ge0},~\{\g_k\}_{k\ge0}\in L^\be(\ell^2)$ and 
$\{\b_{jk}\}_{j,k\ge0}\in L^\be(\B)$, where $\B$ is the space of infinite matrices that induce bounded linear operators on $\ell^2$. We refer the reader to \cite{Pi} for Haagerup tensor products. It is well known (see \cite{JTT}) that for a function $\Psi$ satisfying \rf{Htp}
and for bounded linear operators $T$ and $R$, one can define the triple operator
integral
\bay
\label{toi}
W=\int_{\X_1}\int_{\X_2}\int_{\X_3}\Psi(x_1,x_2,x_3)\,dE_1(x_1)T\,dE_2(x_2)R\,dE_3(x_3)
\ey
and
\bay
\label{Hto}
\|W\|\le\|\Psi\|_{L^\be\otimes_{\rm h}\!L^\be\otimes_{\rm h}\!L^\be}\|T\|\cdot\|R\|,
\ey
where
$$
\|\Psi\|_{L^\be\otimes_{\rm h}\!L^\be\otimes_{\rm h}\!L^\be}\df\inf\left\{
\|\{\a_j\}_{j\ge0}\|_{L^\be(\ell^2)}\|\{\b_{jk}\}_{j,k\ge0}\|_{L^\be(\B)}
\|\{\g_k\}_{k\ge0}\|_{L^\be(\ell^2)}\right\},
$$
the infimum being taken over all representations of $\Psi$ in the form \rf{Htp}.
Note that this extends the definition of triple operator integrals given in \cite{Pe3} for projective tensor products of $L^\be$ spaces.

We would like to define triple operator integrals of the form \rf{toi} in the case
when one of the operators $T$ or $R$ is of trace class and find conditions on $\Psi$ under which the operator $W$ must be in $\bS_1$.

{\bf Definition.}
{\it A function $\Psi$ is said to belong to the tensor product $L^\be(E_1)\otimes_{\rm h}\!L^\be(E_2)\otimes^{\rm h}\!L^\be(E_3)$ if it admits a representation
\bay
\label{yaH}
\Psi(x_1,x_2,x_3)=\sum_{j,k\ge0}\a_j(x_1)\b_{k}(x_2)\g_{jk}(x_3),\quad x_j\in\X_j,
\ey
with $\{\a_j\}_{j\ge0},~\{\b_k\}_{k\ge0}\in L^\be(\ell^2)$ and 
$\{\g_{jk}\}_{j,k\ge0}\in L^\be(\B)$. For a bounded linear operator $R$ and 
for a trace class operator $T$, we define the triple operator integral
$$
W=\int_{\X_1}\int_{\X_2}\int_{\X_3}\Psi(x_1,x_2,x_3)\,dE_1(x_1)T\,dE_2(x_2)R\,dE_3(x_3)
$$
as the following continuous linear functional on the class of compact operators:
\bay
\label{fko}
Q\mapsto
\trace\left(\left(
\int_{\X_1}\int_{\X_2}\int_{\X_3}
\Psi(x_1,x_2,x_3)\,dE_2(x_2)R\,dE_3(x_3)Q\,dE_1(x_1)
\right)T\right).
\ey
}

The fact that the linear functional \rf{fko} is continuous on the class of compact operators is a consequence of inequality \rf{Hto}, which also implies the following estimate:
$$
\|W\|_{\bS_1}\le\|\Psi\|_{L^\be\otimes_{\rm h}\!L^\be\otimes^{\rm h}\!L^\be}
\|T\|_{\bS_1}\|R\|,
$$
where $\|\Psi\|_{L^\be\otimes_{\rm h}\!L^\be\otimes^{\rm h}\!L^\be}$ is the infimum of $\|\{\a_j\}_{j\ge0}\|_{L^\be(\ell^2)}\|\{\b_k\}_{k\ge0}\|_{L^\be(\ell^2)}
\|\{\g_{jk}\}_{j,k\ge0}\|_{L^\be(\B)}$
over all representations in \rf{yaH}.

Similarly, suppose that $\Psi\in L^\be(E_1)\otimes^{\rm h}\!L^\be(E_2)\otimes_{\rm h}\!L^\be(E_3)$, i.e., $\Psi$ admits a representation
$$
\Psi(x_1,x_2,x_3)=\sum_{j,k\ge0}\a_{jk}(x_1)\b_{j}(x_2)\g_k(x_3)
$$
where $\{\b_j\}_{j\ge0},~\{\g_k\}_{k\ge0}\in L^\be(\ell^2)$, 
$\{\a_{jk}\}_{j,k\ge0}\in L^\be(\B)$, $T$ is a bounded linear operator, and $R\in\bS_1$. Then
the continuous linear functional 
$$
Q\mapsto
\trace\left(\left(
\int_{\X_1}\int_{\X_2}\int_{\X_3}
\Psi(x_1,x_2,x_3)\,dE_3(x_3)Q\,dE_1(x_1)T\,dE_2(x_2)
\right)R\right)
$$
on the class of compact operators determines a trace class operator 
$$
W\df\int_{\X_1}\int_{\X_2}\int_{\X_3}\Psi(x_1,x_2,x_3)\,dE_1(x_1)T\,dE_2(x_2)R\,dE_3(x_3).
$$
Moreover,
$$
\|W\|_{\bS_1}\le
\|\Psi\|_{L^\be\otimes^{\rm h}\!L^\be\otimes_{\rm h}\!L^\be}
\|T\|\cdot\|R\|_{\bS_1}.
$$

Note that the above definitions of triple operator integrals extend the definition given in \cite{Pe3} in terms of the projective tensor product of the $L^\be$ spaces.

We would like to obtain a sufficient condition on a function $f$ on $\R^2$ under which 
$$
\|f(A_1,B_1)-f(A_2,B_2)\|_{\bS_1}
\le\const(\|A_1-A_2\|_{\bS_1}+\|B_1-B_2\|_{\bS_1}),
$$
whenever $(A_1,B_1)$ and $(A_2,B_2)$ are pairs of (not necessarily commuting) self-adjoint operators such that $A_1-A_2\in\bS_1$
and $B_1-B_2\in\bS_1$.
Clearly, it suffices to verify the following inequalities: 
$$
\|f(A_1,B)-f(A_2,B)\|\le\const\|A_1-A_2\|_{\bS_1}\quad\mbox{and}\quad
\|f(A,B_1)-f(A,B_2)\|\le\const\|B_1-B_2\|_{\bS_1}.
$$


\begin{thm}
\label{rrp}
Let $f$ be a bounded function on $\R^2$ whose Fourier transform is supported 
in the ball $\{\xi\in\R^2:~\|\xi\|\le1\}$. Then
$$
\frac{f(x_1,y)-f(x_2,y)}{x_1-x_2}=
\sum_{j,k\in\Z}\frac{\sin(x_1-j\pi)}{x_1-j\pi}\cdot\frac{\sin(x_2-k\pi)}{x_2-k\pi}
\cdot\frac{f(j\pi,y)-f(k\pi,y)}{j\pi-k\pi}.
$$
Moreover,
$$
\sum_{j\in\Z}\frac{\sin^2(x_1-j\pi)}{(x_1-j\pi)^2}
=\sum_{k\in\Z}\frac{\sin^2(x_2-k\pi)}{(x_2-k\pi)^2}=1,
\quad x_1~x_2\in\R,
$$
and
$$
\sup_{y\in\R}\left\|\left\{\frac{f(j\pi,y)-f(k\pi,y)}{j\pi-k\pi}
\right\}_{j,k\in\Z}\right\|_\B\le\const\|f\|_{L^\be(\R)}.
$$
\end{thm}

Note that if $j=k$, we assume that $(f(j\pi,y)-f(k\pi,y))(j\pi-k\pi)^{-1}
=\frac{\partial f}{\partial x}(j\pi,y)$.

Theorem \ref{rrp} implies the following result:

\begin{thm}
\label{yadA}
Let $f$ be a bounded function on $\R^2$ whose Fourier transform is supported 
in the ball $\{\xi\in\R^2:~\|\xi\|\le\s\}$. Then the function $\Psi$ defined by
$$
\Psi(x_1,x_2,y)=\frac{f(x_1,y)-f(x_2,y)}{x_1-x_2},\quad x_1,~x_2,~y\in\R,
$$
belongs to the tensor product $L^\be(\R)\otimes_{\rm h}\!L^\be(\R)\otimes^{\rm h}\!L^\be(\R)$,
$$
\|\Psi\|_{L^\be\otimes_{\rm h}\!L^\be\otimes^{\rm h}\!L^\be}
\le\const\s\|f\|_{L^\be(\R^2)}
$$
and
\bay
\label{A1A2}
f(A_1,B)-f(A_2,B)=\int_\R\int_\R\int_\R\frac{f(x_1,y)-f(x_2,y)}{x_1-x_2}
\,dE_{A_1}(x_1)(A_1-A_2)\,dE_{A_2}(x_2)\,dE_B(y),
\ey
whenever $A_1$, $A_2$, and $B$ are self-adjoint operators such that 
$A_1-A_2\in\bS_1$.
\end{thm}

In a similar way we can prove the following fact:

\begin{thm}
\label{yadB}
Under the same hypotheses on $f$, the function $\Psi$ defined by
$$
\Psi(x,y_1,y_2)=\frac{f(x,y_1)-f(x,y_2)}{y_1-y_2},\quad x,~y_1,~y_2\in\R,
$$
belongs to the tensor product $L^\be(\R)\otimes^{\rm h}\!L^\be(\R)\otimes_{\rm h}\!L^\be(\R)$, $
\|\Psi\|_{L^\be\otimes^{\rm h}\!L^\be\otimes_{\rm h}\!L^\be}
\le\const\s\|f\|_{L^\be(\R^2)}$,
and
\bay
\label{B1B2}
f(A,B_1)-f(A,B_2)=\int_\R\int_\R\int_\R\frac{f(x,y_1)-f(x,y_2)}{y_1-y_2}
\,dE_A(x)\,dE_{B_1}(x_1)(B_1-B_2)\,dE_{B_2}(y_2),
\ey
whenever $A$, $B_1$, and $B_2$ are self-adjoint operators such that 
$B_1-B_2\in\bS_1$.
\end{thm}

Theorems \ref{yadA} and \ref{yadB} imply the main result of this section:

\begin{thm}
Let $f\in B_{\be,1}^1(\R^2)$. Then
\bay
\label{liya}
\|f(A_1,B_1)-f(A_2,B_2)\|_{\bS_1}
\le\const\|f\|_{B_{\be,1}^1}(\|A_1-A_2\|_{\bS_1}+\|B_1-B_2\|_{\bS_1}),
\ey
whenever $(A_1,B_1)$ and $(A_2,B_2)$ are pairs of self-adjoint operators, $A_1-A_2\in\bS_1$ and $B_1-B_2\in\bS_1$.
\end{thm}

We have defined functions $f(A,B)$ for $f$ in 
$B_{\be,1}^1(\R^2)$ only for bounded self-adjoint operators $A$ and $B$. 
However, formulae \rf{A1A2} and \rf{B1B2} allow us to define the difference $f(A_1,B_1)-f(A_2,B_2)$ in the case when $f\in B_{\be,1}^1(\R^2)$ and the 
self-adjoint operators $A_1,\,A_2,\,B_1,\,B_2$ are possibly unbounded once we know that the pair $(A_2,B_2)$ is a trace class perturbation of the pair $(A_1,B_1)$.
Moreover, inequality \rf{liya} also holds for such operators.

\

\section{\bf Lipschitz type estimates in the operator norm}

\

The main result of this section shows that unlike in the case of commuting pairs of self-adjoint operators, the condition $f\in B_{\be,1}^1(\R^2)$ does not imply 
Lipschitz type estimates in the operator norm.

\begin{thm}
There is no positive number $M$ such that 
$$
\|f(A_1,B_1)-f(A_2,B_2)\|\le M\|f\|_{L^\be(\R^2)}(\|A_1-A_2\|+\|B_1-B_2\|)
$$ 
for all bounded functions $f$ on $\R^2$ with Fourier transform supported in $[-1,1]^2$ and for all finite rank self-adjoint operators $A_1,\,A_2,\,B_1,\,B_2$.
\end{thm}

\medskip

{\bf Construction.} Let $\{f_j\}_{1\le j\le N}$ and $\{g_k\}_{1\le k\le N}$
be orthonormal systems. Consider the rank one projections $P_j$ and $Q_j$ defined by
$$
P_j u=(u,f_j)f_j\quad\mbox{and}\quad Q_ju=(u,g_j)g_j,\quad 1\le j\le N.
$$
We define the self-adjoint operators $A_1$, $A_2$, and $B=B_1=B_2$ by
$$
A_1=\sum_{j=1}^N4j\pi P_j,\quad A_2=\sum_{j=1}^N(4j+2)\pi P_j,\quad\mbox{and}\quad
B=\sum_{k=1}^N4k\pi Q_k.
$$
Clearly, $\|A_1-A_2\|\le2\pi$.
Put
$$
\f(x,y)=4\cdot\frac{1-\cos x}{x^2}\cdot\frac{1-\cos y}{y^2},\quad x,~y\in\R.
$$
Given a matrix $\{\tau_{jk}\}_{1\le j,k\le N}$, we define 
the function $f$ by
$$
f(x,y)=\sum_{1\le j,k\le N}\tau_{jk}\f(x-4\pi j,y-4\pi k).
$$ 
It is easy to see that 
$$
\|f\|_{L^\be(\R^2)}\le\const\max_{1\le j,k\le N}|\tau_{jk}|,\quad
f(4j\pi,4k\pi)=\tau_{jk},\quad f((4j+2)\pi,4k\pi)=0,\quad1\le j,k\le N,
$$
and the Fourier transform of
$f$ is supported in $[-1,1]\times[-1,1]$.
One can easily verify that
$$
f(A_1,B)=
\sum_{1\le j,k\le N}\tau_{jk}P_jQ_k,\quad\mbox{while}\quad
f(A_2,B)=\0.
$$
It can easily be shown that the supremum of $\|f(A_1,B)\|$ over all $P_j$ and $Q_k$ is the Schur multiplier norm of the matrix $\{\tau_{jk}\}_{1\le j,k\le N}$. This norm can be made arbitrarily large as $N\to\be$, while assumption $|\tau_{jk}|\le1$, $1\le j,k\le N$, can still hold.

%
\medskip

{\bf Remark.} The above construction shows that there exist a function $f$ on $\R^2$ whose Fourier transform is supported in $[-1,1]^2$ such that $\|f\|_{L^\be(\R)}\le1$ and self-adjoint operators of finite rank $A_1$, $A_2$, $B_1$, $B_2$ such that $\|A_1-A_2\|\le2\pi$ 
and $\|B_1-B_2\|\le2\pi$, but $\|f(A_1,B_1)-f(A_2,B_2)\|$ is greater than any given positive number. In particular, the fact that $f$ is a H\"older function of order $\a\in(0,1)$ does not imply the estimate 
$\|f(A_1,B_1)-f(A_2,B_2)\|\le\const(\|A_1-A_2\|^\a+\|B_1-B_2\|^\a)$.


\medskip

We conclude the paper with a theorem that can be deduced from the results of this section.

\begin{thm}
There are spectral measure $E_1$, $E_2$ and $E_3$ on Borel subsets of $\R$, a function $\Psi$ in the Haagerup tensor product $L^\be(E_1)\otimes_{\rm h}\!L^\be(E_2)\otimes_{\rm h}\!L^\be(E_3)$ and an operator $Q$ in $\bS_1$ such that
$$
\iiint\Psi(x_1,x_2,x_2)\,dE_1(x_1)\,dE_2(x_2)Q\,dE_3(x_3)\not\in\bS_1.
$$
\end{thm}

\end{document}

%% file: CRinput.tex
\renewcommand{\a}{\alpha}
\renewcommand{\b}{\beta}
\newcommand{\g}{\gamma}

\newcommand{\s}{\sigma}

\newcommand{\f}{\varphi}

\newcommand{\B}{{\mathcal B}}

\newcommand{\X}{{\mathcal X}}

\newcommand{\R}{{\Bbb R}}
\newcommand{\Z}{{\Bbb Z}}

\newcommand{\0}{{\boldsymbol{0}}}

\newcommand{\bS}{{\boldsymbol S}}

\newcommand{\rf}[1]{(\ref{#1})}

\newcommand{\df}{\stackrel{\mathrm{def}}{=}}

\newcommand{\trace}{\operatorname{trace}}

\newcommand{\const}{\operatorname{const}}

\newcommand{\eeq}{\end{equation}}
\newcommand{\beq}{\begin{equation}}
\newcommand{\bay}{\begin{eqnarray}}
\newcommand{\ba}{\begin{align*}}
\newcommand{\ea}{\end{align*}}
\newcommand{\ey}{\end{eqnarray}}
\newcommand{\bey}{\begin{eqnarray*}}
\newcommand{\eey}{\end{eqnarray*}}

\newcommand{\be}{\infty}

\newtheorem{thm}{\hspace{\parindent}Theorem}[section]